\documentclass[12pt]{amsart}
\usepackage{latexsym, amsmath,amssymb}
\usepackage{hyperref}
\usepackage{xcolor}
\hypersetup{
    colorlinks,
    linkcolor={red!50!black},
    citecolor={blue!50!black},
    urlcolor={blue!80!black}
}

\usepackage{graphicx}

 \numberwithin{equation}{section}

\def\XXint#1#2#3{{\setbox0=\hbox{$#1{#2#3}{%
\int}$ }
\vcenter{\hbox{$#2#3$ }}\kern-.6\wd0}}

\renewcommand{\epsilon}{\varepsilon}
\newtheorem{theorem}{Theorem}

\newtheorem{lemma}[theorem]{Lemma}
\newtheorem{corr}[theorem]{Corollary}

\newtheorem{proposition}[theorem]{Proposition}
\newtheorem{deff}[theorem]{Definition}

\newcommand{\bth}{\begin{theorem}}
\newcommand{\ble}{\begin{lemma}}
\newcommand{\bcor}{\begin{corr}}

\newcommand{\bdeff}{\begin{deff}}

\newcommand{\bprop}{\begin{proposition}}
\newcommand{\ele}{\end{lemma}}
\newcommand{\ecor}{\end{corr}}
\newcommand{\edeff}{\end{deff}}

\numberwithin{theorem}{section}

\newcommand{\eprop}{\end{proposition}}

\renewcommand{\Pi}{\varPi}

\renewcommand{\epsilon}{\varepsilon}

\begin{document}

\title[An abstract generalization of K\"all\'en-Kr\"oner's iteration]
{An abstract generalization\\ of K\"all\'en-Kr\"oner's iteration}
\author[A. D. Mart\'inez]{\'Angel D. Mart\'inez}
\address{Institute for Advanced Study, Fuld Hall 412, 1 Einstein Drive, Princeton, NJ 08540, United States of America}
\email{amartinez@ias.edu}

\begin{abstract}
We devote this paper to provide an abstract generalization of an iteration originally due to K\"all\'en, and revisited later by Kr\"oner, that might be of independent interest. An application to prove a Nash-Kuiper theorem for the existence of $C^{1,1-\epsilon}$ isometric embeddings in codimension $\frac{1}{2}n(n+1)$ can be found in \cite{paper3.21}.
\end{abstract}

\maketitle

\section{\textbf{Introduction}}

We will present an iteration originally present in the work of K\"all\'en that provides an useful quantitative almost fixed point theorem in conection with its original application to the Borisov-Gromov problem about the existence (or not) of $C^{1,\alpha}$ isometric embeddings of compact Riemannian manifolds. Both K\"all\'en and Kr\"oner's used the iteration implicitly which might make it slightly cumbersome to handle or conceptualize. The main result in this paper isolates a sensible abstract general set of hypothesis under which the iteration works. We hope this will help the iteration to reach other potential uses in other applications as well. 

We refrain to present Nash's iteration in detail here and refer the reader to the delightful original paper \cite{N} (cf. \cite{CLS}). Let us summarize it by saying that at each step in Nash's iteration one needs to produce a highly oscillatory perturbation that approximates a (metric) tensor $T$ by a tensor of the form $b(a,a)+r(a)$. In this case $a$ denotes a vector valued function (containing the coefficients of the perturbation in a specific basis), $b$ is bilinear in the $a$'s with tensor valued functions as codomain and $r$ is a function from vectors valued functions to tensor valued functions that depends on the variable $a$ (actually, in the application, it is bilinear too but not completely in $a$ since it might contain linear terms too). One then tries to minimize the tensor error $E=T-b(a,a)-r(a)$ in some suitable norm.  This kind of procedure is at the heart of the Nash-Kuiper iteration where a small error is commited after adding an oscillatory perturbation at frequency $\lambda$ to a mollification of the previous step at scale $\ell$. At the end of the day the smaller the error the smoother the embedding will be. In the best of the worlds one would employ a fixed point argument to find the best choice for $a$. Unfortunately, there is a loss of derivatives inherent to the equation at hand. Quite surprisingly, the K\"all\'en-Kr\"oner's iteration improves the error term quantitatively by worsening some of the involved constants.  The upshot is that the constants will be fixed in the Nash-Kuiper iteration and can be forgotten right away.

Let us fix some notation now. For the time being, and for the sake of concreteness in the exposition, let $b$ be a fixed bilinear function (from pairs of vectors to tensors) and $r$ another one that can be decomposed as a linear combination 
\[r=r_1+r_2+r_3+r_4\]
with $r_1$ linear in $a$, and rest are bilinear in $a$, but maybe depending on its derivatives in a suitable way. In particular, notice that $r(0)=0$. More concretely, we will require each $r_i$ to satisfy certain mild conditions. For the linear term we impose the following type of bounds: the linear term satisfies
\[\|r_1(a)\|_k\leq  \frac{C_k}{\lambda \ell}\sum_{j=0}^k\|a\|_j\lambda^{k-j}.\]
Similarly, we allow quadratic terms satisfying
\[\|r_2(a, b)\|_k\leq \frac{C_k}{\lambda^{2}\ell^{2}}\sum_{j_1+j_2=0}^k\|a\|_{j_1}\|b\|_{j_2}\lambda^{k-j_1-j_2},\]
\[\|r_3(\nabla a,\nabla b)\|_k\leq\frac{C_k}{\lambda^{2}}\sum_{j_1+j_2=0}^k\|\nabla a\|_{j_1}\|\nabla b\|_{j_2}\lambda^{k-j_1-j_2}\]
or
\[\|r_4(a, b)\|_k\leq \frac{C_k}{\lambda^{2}\ell}\sum_{j_1+j_2=0}^k\|a\|_{j_1+1}\|b\|_{j_2}\lambda^{k-j_1-j_2}.\]

In the application to the isometric embedding  it is showed that one can essentially gain the $\ell$ factor from the mollification. The reader can think of this $\ell$ as a parameter that might or might not be directly related to the mollification parameter in future applications.

This estimates will need to be checked in the application at hand which is essentially done in Kr\"oner's Master thesis for the isometric embedding. He actually splits the perturbation into five possible forms and bound them separately. This is tedious and might obscure the essence of the iteration which we are trying to highlight here, probably paying the price of obscuring it in the abstraction (cf. \cite{K}). We believe though that this abstract version might help the community applying this simple and powerful lemma to other problems. The iteration's first step is a {\em deus ex machina} choice for $a$ such that $T=b(a,a)$. Notice this does not take $r$ into account. This is already in Nash's original work which bounds the error terms as $O(\lambda^{-1})$. As a consequence $r$ becomes a rather small error term choosing $\lambda$ large enough at each step. This works wonderfully if one is interested in $C^1$ or slightly better embeddings but does not provide sharp results. Let us mention that a variant of the method due to K\"allen takes part of this term into account at the expense of working in high codimension (cf. \cite{Ka} where the global case is considered or  \cite{K} for the local case). The method, nevertheless, does not deal with self interaction terms satisfying
\[\|r_5(a, b)\|_k\leq \frac{C_k}{\lambda\ell}\sum_{j_1+j_2=0}^k\|a\|_{j_1+1}\|b\|_{j_2}\lambda^{k-j_1-j_2}\]
and, as a consequence, the K\"all\'en-Kr\"oner iteration still needs high codimension to make errors of this type vanish {\em a priori} from purely geometrical reasons. This last  term must be compared with the linear type term $r_1$ for which a $(\lambda\ell)^{-1}$ gain is enough.

We still need to fix some more notation let $F$ be a function from tensors to vectors that asigns to each $T$ a solution $a$ to the equation $T=b(a,a)$, i.e. $F(T)=a$. In practical situations we can ensure the existence of such an inverse, smoothly, in the neighborhood of a fixed tensor $T_0$. In some applications this will be the identity and the neighbourhood is in the $C^0$ topology. In the hypothesis for the K\"all\'en-Kr\"oner's iteration  we will  need to control more derivatives to ensure that also derivatives of $F$ evaluated at any tensor $T$ close to $T_0$ are comparable with values at $T_0$ which is fixed. As the reader can check in the hypothesis of Lemma \ref{inductivekroner} below we will assume that we are working in such a scenario from now on. 

The statement of K\"all\'en-Kr\"oner is rather long and technical, we postpone it to the last section (cf. Theorem \ref{kr}).

In the next section we provide a concrete statement of the inductive step. We present its proof in Section \ref{proof}. In the last section we state and prove the main result in a rather stronger (and more abstract) form which adapts better to its application to the Borisov-Gromov problem.

\section{Statement of the inductive step}

Let us state first the following inductive procedure (which is an abstract generalization of K\"all\'en's ingenious idea).

\begin{lemma}\label{inductivekroner}
Let $F$ be the inverse of $b$ defined in a $C_F^{-1}$-neighbourhood of $T_0$ be such that its first $k_0+1$ derivatives satisfy
\[\|F(T)\|_k\leq C_F\left(\|T\|_k+\frac{\lambda^k}{\lambda\ell}\right)\]
and
\[\|F(T)-F(T')\|_k\leq C_{F}(\|T-T'\|_k+(\|T'\|_k+\lambda^k)\|T-T'\|_0)\]
for tensors $T$ and $T'$ in a $3C_F^{-1}$-neighbourhood of $T_0$  for some fixed $k_0\in\mathbb{N}$, $b$ and $r^{i-1}$, $r^{i}$ and $r^{i+1}$ as already described in the introduction. Suppose that for two fixed constants $\lambda,\ell>0$ such that $\lambda\ell$ is large enough (depending on $C_F$). 

Given another fixed tensor $T$ in the $(3C_F)^{-1}$-neighbourhood of $T_0$ satisfying 
\[\|T\|_k\leq C\frac{\lambda^k}{\lambda\ell}\]
for $k\geq 1$ and a vector $a^{(i)}=F(T-r^{i-1}(a^{(i-1)}))$ providing an approximation  error 
\[E_i=T-b(a^{(i)},a^{(i)})-r^{i}(a^{(i)})\]
and satisfying
\begin{itemize}
\item[(1)] $\|a^{(i)}\|_0\leq C$.
\item[(2)] For any $1\leq k\leq k_0-i$ the norm 
\[\|a^{(i)}\|_k\leq C\frac{\lambda^k}{\lambda\ell},\]
\item[(3)] while the error satisfies
\[\|E_i\|_k\leq C\frac{\lambda^k}{(\lambda\ell)^i},\]
\item[(4)] and, finally, we also have the bounds
\[\|r^{i}(a^{(i)})\|_k\leq C\frac{\lambda^k}{\lambda\ell}\]
and
\[\|r^{i-1}(a^{(i-1)})\|_k\leq C\frac{\lambda^k}{\lambda\ell}.\]
\item[(5)] Furthermore, suppose inductively that, for every step $i$, the $r^i$ satisfy the bounds already introduced and
\[\|r_1^{i+1}(a)-r_1^{i}(a)\|_k\leq  \frac{C_k}{\lambda \ell}\sum_{j=0}^k\|a\|_j\lambda^{k-j},\]
\[\|r_2^{i+1}(a,b)-r_2^{i}(a,b)\|_k\leq \frac{C_k}{\lambda^{2}\ell^{2}}\sum_{j_1+j_2=0}^k\|a\|_{j_1}\|b\|_{j_2}\lambda^{k-j_1-j_2},\]
\[\|r_3^{i+1}(\nabla a,\nabla b)-r_3^{i}(\nabla a,\nabla b)\|_k\leq \frac{C_k}{\lambda^{2}}\sum_{j_1+j_2=0}^k\|\nabla a\|_{j_1}\|\nabla b\|_{j_2}\lambda^{k-j_1-j_2},\]
\[\|r_4^{i+1}(a,b)-r_4^{i}(a,b)\|_k\leq \frac{C_k}{\lambda^{2}\ell}\sum_{j_1+j_2=0}^k\|a\|_{j_1+1}\|b\|_{j_2}\lambda^{k-j_1-j_2},\]
and $\|r^{i}(a)\|_0\leq C_0$ uniformly in $i$ for $\|a\|_0\leq C$.
\end{itemize}
Then, if we define $a^{(i+1)}:=F(T-r^{i}(a^{(i)}))$, it satisfies the above estimates (1)-(4) with $i$ replaced by $i+1$.
\end{lemma}

\textsc{Remark:} the constants involved implicitly depend on $i$, $k_0$ and $C_F$ and, as is customary, might change from line to line. The first bound in (4) is redundant as it follows from (1), (2) and the hypothesis on $r$. The hypothesis on $F$ corresponds to H4 (ii) in \cite{Ka} where an extra condition $\|T'\|_0\leq C'$ for some fixed constant follows from our hypothesis with $C'=C'(C_F,T_0)$. The relevance of this seemingly vacous hypothesis is that it allows to change $F$ from step to step which will be useful in some applications (cf. Theorem \ref{kr}). 

In the case of the Borisov-Gromov problem, as we have already mentioned, $T$ will be the target metric defect that we want to reduce, $\lambda$ corresponds to the next frequency $\lambda_{q+1}$ in the iteration process and $\ell$ would correspond to the mollification parameter  at the begining of each stage. Notice that the hypothesis $\lambda\ell$ big enough is sensible since by (3) that improves the error term (which is our objective) while the bounds on $a$ remain virtually the same (as they should). The condition on $T$ is also sensible since in applications it will be a mollification of a fixed tensor for which a much better estimate would be available.

Notice that in general there is a loss of derivatives inherent to the lemma. Indeed, in general $r(a^{(i)})$ contains derivatives of $a^{(i)}$, say one derivative, then to control $k$ derivatives of $a^{(i+1)}=F(T-r(a^{(i)}))$ one needs to control $k+1$ derivatives of $a^{(i)}$. This produces a cascade of a loss of derivatives in the induction if we assume we are controlling only finitely many from the very begining of the process as in (2) above.

We must warn the reader that (5) is a hypothesis that must be checked by some different means. In fact, K\"all\'en-Kr\"oner's work deal with error terms that do not change from step to step for which this condition is empty. This version is not just an abstract rephrasing of his idea but a strenghtening too. The applications will appear elsewhere.

\section{Proof of Lemma \ref{inductivekroner}}\label{proof}

Implicitly we are assuming such $a^{(i+1)}$ are well defined. We begin quickly checking that they are indeed well defined. To do so we need first to check that $T-r^{i}(a^{(i)})$ is close to $T_0$. Since $T$ is already close to $T_0$ this will happen provided $r^{i}(a^{(i)})$ is small enough to be in the neighbourhood of radius, say, $1/C_F$, where $F$ is defined. In fact, this is a consequence of $\|T-T_0\|_0\leq 1/(3C_F)$ and $\|r^{i}(a^{(i)})\|_0\leq 1/(3C_F)$, the last of which follows by imposing $\lambda\ell$ to be large enough (depending on $C_F$).

The bounds on (1) and (2) follow from the definition of $a^{(i+1)}$ through $F$ using (1), (2), (4) for $a^{(i)}$ and the hypothesis on the tensor $T$ as follows. Indeed,
\[\|a^{(i+1)}\|_0=\|F(T-r^{i}(a^{(i)}))\|_0\]
which should be similar, by continuity of $F$, to the size of $F(T_0)$ and therefore uniformily bounded in its $C_F^{-1}$-neighbourhood. On the other hand using the hypothesis on $F$ one gets
\[\begin{aligned}
\|a^{(i+1)}\|_k&=\|F(T-r^{i}(a^{(i)}))\|_k\\
&\leq C_{F}\left(\|T-r^i(a^{(i)})\|_k+\frac{\lambda^k}{\lambda\ell}\right)\\
&\leq C_F\left(\|T\|_k+\|r^i(a^{(i)})\|_k+\frac{\lambda^k}{\lambda\ell}\right).
\end{aligned}\]
Now it is clear that we can use the inductive hypothesis to close this part of the argument.

To show (3) let us observe first that by definition
\[E_{i+1}=T-b(a^{(i+1)},a^{(i+1)})-r^{i+1}(a^{(i+1)})\]
and also by definition of $a^{(i+1)}$ we know that it solves
\[T-r^{i}(a^{(i)})=b(a^{(i+1)},a^{(i+1)}).\]
Plain substitution gives the remarkable identity
\begin{equation}\label{remarkable}
E_{i+1}=r^{i}(a^{(i)})-r^{i+1}(a^{(i+1)}).
\end{equation}
We will use this to prove (3) and (4). Before doing so recall that by hypothesis the $r$ are combinations of functions either linear in $a$ or bilinear in $a$ (and maybe some of its derivatives). The elementary identity $(a^2-b^2=a(a-b)-b(b-a)$ suggests we study the norms of differences, namely
\[a^{(i+1)}-a^{(i)}=F(T-r^{i}(a^{(i)}))-F(T-r^{i-1}(a^{(i-1)})),\]
as an intermediate step. The identity follows by definition. Using the second hypothesis on $F$ we can bound
\[\|a^{(i+1)}-a^{(i)}\|_k=\|F(T-r^{i}(a^{(i)}))-F(T-r^{i-1}(a^{(i-1)}))\|_k\]
from above by
\[ C_{F}(\|r^{i}(a^{(i)})-r^{i-1}(a^{(i-1)})\|_{k}+\|T-r^{i-1}(a^{(i-1)})\|_k\|r^{i}(a^{(i)})-r^{i-1}(a^{(i-1)})\|_0).\]
Using the identity \ref{remarkable} for $i$ and the hypothesis on $F$ one gets
\[\|a^{(i+1)}-a^{(i)}\|_k\leq C_{F}(\|E_i\|_k+(\|T\|_k+\|r^{i-1}(a^{(i-1)})\|_k+\lambda^k)\|E_i\|_0)\]
which are quantities we control from the induction hypothesis (3) and (4), respectively. Using them we finally get
\[C_F\left(C\frac{\lambda^k}{(\lambda\ell)^i}+\left(C\frac{\lambda^k}{(\lambda\ell)}+C\frac{\lambda^k}{(\lambda\ell)}+\lambda^k\right)C\frac{1}{(\lambda\ell)^i}\right).\]
Summarizing we get
\begin{equation}\label{claim}
\|a^{(i+1)}-a^{(i)}\|_k\leq C_{k,F}\frac{\lambda^k}{(\lambda\ell)^i}.
\end{equation}

Now we will concentrate on understanding what happens with the norm of $E_{i+1}$. Taking advantage of the identity \ref{remarkable} one can split this into four parts each of which is a difference of $r_i$. The linear part for example can be bounded as
\[\|r^{i}_1(a^{(i)})-r^{i+1}_1(a^{(i+1)})\|_k\leq \frac{C_k}{\lambda\ell}\sum_{j=0}^k\|a^{(i)}-a^{(i+1)}\|_j\lambda^{k-j}\leq\frac{C_k\lambda^k}{(\lambda\ell)^{i+1}}\]
which follows immediately as a consequence of the elementary
\[r^{i}_1(a^{(i)})-r^{i}_1(a^{(i+1)})+r^{i}_1(a^{(i+1)})-r^{i+1}_1(a^{(i+1)})\]
together with the extra mild hypothesis and  equation \ref{claim}. Similarly for the bilinear part of type $r_2$ let us observe that
\[r^{i+1}_2(a^{(i+1)},a^{(i+1)})-r^{i}_2(a^{(i)},a^{(i)})\]
equals
\[r^{i}_2(a^{(i+1)},a^{(i+1)}-a^{(i)})-r^{i}_2(a^{(i)}-a^{(i+1)},a^{(i)})\]
plus the correction
\[r^{i+1}_2(a^{(i+1)},a^{(i+1)})-r^{i}_2(a^{(i+1)},a^{(i+1)}).\]
The mild condition on $r^{i}_2$ allows to bound the first using equation \ref{claim} to control the differences  and (2) from the hypothesis for $i$ and $i+1$ (which is already proved inductively). The second can be handled using the hypothesis directly to conclude that 
\[\|r^{i+1}_2(a^{(i+1)},a^{(i+1)})-r^{i}_2(a^{(i+1)},a^{(i+1)})\|_k\leq C\frac{\lambda^k}{(\lambda\ell)^{i+1}}\]
The exact same reasoning applies to $r_3$ that depends bilinearly on derivatives of $a$ for which
\[\|r^{i+1}_3(\nabla a^{(i+1)},\nabla a^{(i+1)})-r^{i}_3(\nabla a^{(i)},\nabla a^{(i)})\|_k\leq C\frac{\lambda^{k}}{(\lambda\ell)^{i+1}}\]
holds. We leave to the reader the details for $r^{i}_4$ and comparing them to $r^{i}_5$ which has an unavoidable loss of $\lambda$. This certainly concludes the proof of (3) at stage $i+1$. 

We are only left to prove (4). This can be done bounding $r^{i}$ and $r^{i+1}$ directly using  the bounds (1) and (2), which we already know hold for $a^{(i+1)}$, and the mild hypothesis for these. 

Quite interestingly, we do not need to invoke the specific form of those if we are willing to rely on a complete induction scheme. Indeed, we can always express it using equation \ref{remarkable}
\[r^{i+1}(a^{(i+1)})=r^1(a^{(1)})-\sum_{j=1}^{i+1}E_i\]
and adding the telescoping series that results from the bounds we already have on these. This observation accounts in part for the lack of a condition of type (4) in the statement of Proposition \ref{kr} below. Notice also that to prove the estimate \ref{claim} we do not need (1) nor (2). As a consequence we could have used complete induction and a telescoping argument to show (1) and (2) as a consequence of this.

\textsc{Remark:} notice that the third term might be generalized straightforwardly to higher number of derivatives, i.e.
\[\|r_6(\nabla^sa,\nabla^t b)\|_k\leq C \frac{1}{\lambda^{d}}\sum_{j_1+j_2=0}^k\|\nabla^sa\|_{j_1}\|\nabla^tb\|_{j_2}\lambda^{k-j_1-j_2}\]
where $d=s+t$ denotes de number of derivatives involved  and $s,t\geq 1$. The proof is identical and we omit it. Let us observe though that the lemma would require a small modification this time due to a higher loss of derivatives (depending on an upper bound for the $s$ and $t$).

\section{Statement and proof}\label{statement}

We can finally state

\begin{proposition}[K\"all\'en-Kr\"oner's iteration]\label{kr}
Fix $k_1\in\mathbb{N}$ and $i$, also $b^i$ (with inverse $F^i$) and $r^i$ maps as described above. Then there exist a $k_0=k_0(i,k_1)\in\mathbb{N}$ such that the following is true. Suppose that for two fixed constants $\lambda,\ell>0$ the first $k_0+1$ derivatives of $F^i$  inductively satisfy
\[\|F^i(T)\|_k\leq C_*\left(\|T\|_k+\frac{\lambda^k}{\lambda\ell}\right)\]
and
\[\|F^{i+1}(T)-F^i(T')\|_k\leq C_{*}\left(\|T-T'\|_k+\frac{\lambda^k}{(\lambda\ell)^i}+(\|T'\|_k+\lambda^k)\left(\|T-T'\|_0+\frac{1}{(\lambda\ell)^i}\right)\right)\]
for any $i$ and tensors $T$ and $T'$ in a $3C_*^{-1}$-neighbourhood of $T_0$ and, furthermore, $\lambda\ell$ is large enough (depeding on $C_*$).

Given any tensor $T$ in the $(3C_*)^{-1}$-neighbourhood of $T_0$ satisfying 
\[\|T\|_k\leq C\frac{\lambda^k}{\lambda\ell}\]
for $k\geq 1$. Let $a^{(0)}=0$ and $a^{(1)}=F^1(T)$
\[E_i=T-b^i(a^{(i)},a^{(i)})-r^i(a^{(i)})\]
satisfying
\begin{itemize}
\item[(i)] $\|a^{(1)}\|_0\leq C$.
\item[(ii)] For any $1\leq k\leq k_0-1$ the norm 
\[\|a^{(1)}\|_k\leq C\frac{\lambda^k}{\lambda\ell}\]
\item[(iii)] while the error satisfies
\[\|E_1\|_k\leq C\frac{\lambda^k}{\lambda\ell}.\]
\end{itemize}
Then,  there exist $a^{(i+1)}$ and $r^{i+1}$ satisfying
\begin{itemize}
\item[(1)] $\|a^{(i+1)}\|_0\leq C$.
\item[(2)] For any $1\leq k\leq k_1$ the norm 
\[\|a^{(i+1)}\|_k\leq C'\frac{\lambda^k}{\lambda\ell}\]
\item[(3)] while the error satisfies
\[\|E_{i+1}\|_k\leq C'\frac{\lambda^k}{(\lambda\ell)^{i}}.\]
\item[(4)] and
\[\|a^{(i+1)}-a^{(i)}\|_k\leq C_{k,*}\frac{\lambda^k}{(\lambda\ell)^i}\textrm{ for any $i\geq 1$},\]
\end{itemize}
provided that, by induction, the following inequalities hold 
\[\|r^{i}_1(a^{(i)})-r^{i}_1(a^{(i+1)})\|_k+\|r^{i}_1(a^{(i+1)})-r^{i+1}_1(a^{(i+1)})\|_k\leq \frac{C_k\lambda^k}{(\lambda\ell)^{i+1}},\]
for the linear part, and
\[\|r^{i+1}(a^{(i+1)},a^{(i+1)})-r^{i}(a^{(i+1)},a^{(i+1)})\|\leq \frac{C_k\lambda^k}{(\lambda\ell)^{i+1}}\]
and
\[\|r^{i}(a^{(i+1)},a^{(i+1)}-a^{(i)})\|_k+\|r^{i}(a^{(i)}-a^{(i+1)},a^{(i)})\|_k\leq \frac{C_k\lambda^k}{(\lambda\ell)^{i+1}}\]
for the bilinear part.
\end{proposition}

\textsc{Remark:} notice that the $r^i$ in this proposition can depend on the previous steps. This is convenient in certain situations (cf. \cite{paper3.21}). The original iteration of K\"all\'en-Kr\"oner considers the case $r^i=r$ and $F^i=F$ constant for which the last hypothesis (5) in Lemma \ref{inductivekroner}  are empty. Notice that the true nature of $b^i$ and $r^i$ is somewhat irrelevant to the proof and we are not assuming anything on them (beyond the estimates on its counterpart $F^i$ and $r^i$ respectively) in this last version. The reader can check that the estimates on $F$ needed to prove Lemma \ref{inductivekroner} in this variable case, for instance, equation \ref{claim}, are precisely the ones we state. We should also remark that the iteration (as stated) is not a direct application of the Lemma \ref{inductivekroner}. Indeed, condition (5) has been replaced here by a more general condition  which can be checked to be enough, by inspection of the proof. 

\textsc{Proof:} the hypothesis provide $T=b^1(a^{(1)},a^{(1)})$ exactly and therefore $E_1$ equals $r^{1}(a^{(1)})$. This together with the observation that $r^0(a^{(0)})=0$ explains the absence of the hypothesis (4) as in Lemma \ref{inductivekroner}. Indeed, (iii) immediately implies it. The constant $C'$ depends on $i$, $k_1$, $C_*$ and the constants involved in the conditions imposed on the linear and bilinear parts of $r$. We define $a^{(i+1)}=F^{i+1}(T-r^i(a^{(i)}))$. This together with the fact that the hypothesis (i) and (ii) are the same as (1) and (2) in Lemma \ref{inductivekroner} provide the initial step $i=1$. A complete induction argument produces $a^{(i)}$ provided we take $k_0(i)$ large enough so that $i$ iterations of it still leave us control of the first $k_1$ estimates. The constants then are uniformly bounded for fixed $i$ and therefore $C'$ depends only on $C_*$ and the fixed constants $i$ and $k_1$. This provides (1)-(3). Inspection of the proof shows that (4) follows as a byproduct (cf. equation \ref{claim}).

\section{\textbf{Acknowledgments}}
The author would like express his gratitude to C. De Lellis for drawing his attention to Kr\"oner's master thesis, his encouragement and enlightening discussions around the isometric embedding problem. We are also grateful to L. Sz\'ekelyhidi, Jr. for kindly indicating that the main idea was already present in K\"all\'en's original work. 

This material is based upon work supported by the National Science Foundation under Grant No. DMS-1638352.

\end{document}